\documentclass[a4paper,14pt]{article}
\usepackage[cp1251]{inputenc}
\usepackage[english]{babel}
\usepackage{amsfonts,amssymb,mathrsfs,amscd}
\usepackage[tbtags]{amsmath}
\usepackage{amsfonts,amssymb,mathrsfs,amscd}
\usepackage{graphics}
\usepackage{wrapfig}
\usepackage{euscript}
\usepackage[dvips]{graphicx}

\date{}
\newtheorem{theorem}{Theorem}[section]
\newtheorem{lemma}{Lemma}[section]
\newtheorem{definition}{Definition}[section]

\begin{document}
\renewcommand{\baselinestretch}{0.96}
\renewcommand{\thesection}{\arabic{section}}
\renewcommand{\theequation}{\thesection.\arabic{equation}}
\csname @addtoreset\endcsname{equation}{section}
\Large

\Large
 \centerline{\bf p-adic tight wavelet frames }
  \centerline{\bf S.\,F.~Lukomskii, A.M.Vodolazov}
\vskip 0.5cm
\noindent\large
 N.G.\ Chernyshevskii Saratov State University\\
 LukomskiiSF@info.sgu.ru\\
 vam21@yandex.ru\\

\begin{abstract}
  We propose a simple method to construct   step mask and
  corresponding   step wavelet functions that generate
  tight wavelet frames on the  field $ \mathbb Q_p$ of $p$-adic  numbers.
  To  construct this mask we use three new ideas. First, we  consider only an additive group of the field $ \mathbb Q_p$.
    Second, to construct the mask we   use a special tree. Thirdly, we do not use the principle of unitary extension, we use Pontryagin's principle of duality.
  \\
   Bibliography: 32 titles.
\end{abstract}
\noindent
 keywords: local fields, refinable equation, tight wavelet frames, p-adic numbers, trees.\\
 MSC:Primary 42C40; Secondary 43A75,43A40\\
\Large

\section*{Introduction}

MRA based wavelet tight frames can be viewed as generalization  of the MRA-based orthonormal wavelet. Since the publication of \cite{MS,DI} especially after compactly supported MRA based orthonormal wavelets beind constructed  by \cite{DI},  wavelet analysis   and its applications  have been  one of the most active  research  areas  in applied mathematics.
Algorithms and methods based on wavelet analysis have become powerfulltools in image  signal compression and analysis. Since in zero-dimensional groups the Fourier transform of a step function is again a step function, an interest arose in wavelets on Vilenkin groups  \cite{LBK,LSBG,YuF2}, local fields of positive characteristics \cite{JLJ,BJ1,BJ2,BJ3,BJ4,LVd2,LVB,SSB,FMS,SD}, field of p-adic number \cite{SK,AKSh,KhSh2,ShSk,E1,ES,AES}, and zero-dimensional groups \cite{BRL,LSF}.

The publication of the unitary extension principle \cite{RSh97} in 1997 initiated  a new wave of theoretical development,  as well as exploration of new applications of  MRA-based  tight wavelet  frames.
In the article \cite{FLS} Yu.Farkov, E.Lebedeva, M.Skopina, using the principle of unitary expansion, gave an explicit description of all Vilenkin polynomials generating tight wavelet frames in Vilenkin group.
 In the paper \cite{SD} F.A.Shah, L.Debnaht, suggested some algorithms for constructing tight wavelet frames on local fields of positive characteristic. To construct these frames  authors used the unitary extension principle.
 In the article \cite{ABS} Owais Ahmad, Mohammad Younus Bhat, and Neyaz Ahmad Sheikh, presented   a way to find masks $m_\ell$ that generate fundamental tight frame wavelet.

The notion of p-adic MRA was introduced and a general scheme
for its construction was described in \cite{ShSk}. The first orthogonal wavelet basis was con\-structed  by S. Kozyrev \cite{SK}. A. Yu. Khrennikov, V. M. Shelkovich, M.A.Skopina \cite{KhShSk} described a wide class of orthogonal scaling functions generating  p-adic MRA.

   S.Albeverio,  S.Evdokimov and M.Skopina \cite{AES} proved:
  1)there exists a unique p-adic MRA generated by an orthogonal scaling test function,
  2) there exists non\-orthogonal p-adic MRA
 3) there exists p-adic wavelet frames.  Let us describe the last result in more detail. Let $W_0$ be a wavelet space and there exists functions $
 \psi^{(\nu)}\in L_2({\mathbb Q_p)}$ such that
 $$
 W_0=\overline{{\rm span}\{\psi^{(\nu)}(x\dot{-}h), \nu=1,...,r, h\in H_0\}}.
 $$
   It $\psi^{(\nu)}, \nu=1,...,r$ is a set of compactly  supported wavelet functions for an p-adic MRA $(V_n)$  then the wavelet system
   $$
   \{p^{n/2}\psi^{(\nu)}(\mathcal{A}^nx\dot{-}h), \nu=1,...,r, h\in H_0\}
   $$
   is a frame in $L_2(\mathbb Q_p)$.
   However, there are no methods for constructing nonortho\-gonal  tight p-adic wavelet frames.  \\
   In \cite{LSF}  provides methods for constructing orthogonal wavelets in an arbitrary zero-dimensional group, however, there are no methods for constructing non\-orthogonal tight wavelet  frames.
Recall that a system  $(f_n)_{n=1}^\infty$ is called a tight frame in a Hilbert space $ H$ if for any $f\in H$
$$
\sum_{n\in \mathbb N}|<f,f_n>|^2=\|f\|^2\Leftrightarrow
\sum_{n\in \mathbb N}<f,f_n>f_n=f.
$$
 In this article we  obtain a method for constructing  tight wavelet frames in the field $\mathbb Q_p$  without the principle of unitary extension. In fact, we will discuss this problem in the locally compact zero-dimensional group, which is the additive group of the field $\mathbb Q_p$.

 The paper is organized as follows.

 In Sec. 1, we give the necessary information about locally compact zero-dimensional groups  and prove some auxiliary results.

 In Sec. 2, we discuss the problem of constructing a step refinable  function in the group of p-adic numbers.

  In the third section, we indicate a method for constructing step tight wavelet frames in the group of p-adic numbers and give examples.

\section{MRA and refinable functions in locally compact zero-dimensional Abelian
groups   }
Let $(G,\dot + )$~be a~locally compact zero-dimensional Abelian
group with the topology generated by a~countable system of open
subgroups
$$
\cdots\supset G_{-n}\supset\cdots\supset G_{-1}\supset G_0\supset
G_1\supset\cdots\supset G_n\supset\cdots
$$
where
$$
\bigcup_{n=-\infty}^{+\infty}G_n= G,\quad
 \quad \bigcap_{n=-\infty}^{+\infty}G_n=\{0\},
$$
 $p$ is an order of quotient groups $G_n/G_{n+1}$ for all $n\in\mathbb Z$.
 We will always
assume that $p$ is a prime number. We will name such chain
as \it basic chain. \rm  In this case, a~base of the topology is
formed by all possible cosets~$G_n\dot + g$, $g\in G$.
Given $n\in\mathbb Z$, consider an element $g_n\in G_n\setminus
G_{n+1}$ and fix~it. Then any $g\in G$ has a~unique representation
in the form
\begin{equation}
\label{eq1.01} g=\sum_{n=-\infty}^{+\infty}a_ng_n, \qquad
a_n=\overline{0,p-1}.
\end{equation}
The sum~\ref{eq1.01} contains finite number of terms with negative
subscripts, that~is,
\begin{equation}
\label{eq1.02} g=\sum_{n=m}^{+\infty}a_ng_n, \qquad
a_n=\overline{0,p-1}, \quad a_m\ne 0.
\end{equation}
We will name the system $(g_n)_{n\in \mathbb Z}$ as {\it a basic
system}. The mapping $\lambda:G\to[0,+\infty)$ defined by the equality
$$
\lambda(g)=\sum_{n=m}^{+\infty}a_np^{-n-1}
$$
is called Monna mapping \cite{AKS}.  It is evident that $\lambda(G_n)=[0,p^{-n}]$.
 Therefore, geometrically, the subgroups $ G_n $ can be represented on the number line as follows

   \unitlength=0.80mm
 \begin{picture}(200,45)
 \put(0,20){\line(1,0){180}}
 \qbezier(0,20)(10,25)(20,20)
 \put(40,18){\line(0,1){4}}
 \put(60,18){\line(0,1){4}}
 \put(120,18){\line(0,1){4}}
 \qbezier(120,20)(150,10)(180,20)
 \qbezier(0,20)(30,10)(60,20)
 \put(15,9){\small ${G_{n+1}}$}
 \put(4,26){\small ${G_{n+2}}$}
 \put(55,12){\small $\frac{1}{p^{n+1}} $}
 \put(130,8){\small $ G_{n+1}\dot+(p-1)g_n$}
 \put(90,22){\small $g_n $}
 \put(90,20){\circle*{1} }
 \qbezier(0,20)(90,45)(180,20)
 \put(100,35){\small $G_n $}
 \end{picture}

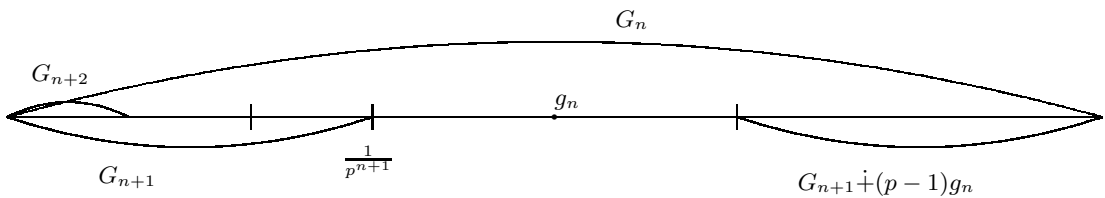
\hskip2cm Figure 1.  Representation of the group $G$ on the number line.
\vskip0.4cm
 Vilenkin groups and groups of $p$-adic numbers (see~\cite[Ch.~1, \S\,2]{AVDR}) are classical examples of zero-dimensional groups .
By $X$ we denote  the collection of the characters of a~group $(G,\dot+ )$; it is
a~group with respect to multiplication, too. Also let
$G_n^\bot=\{\chi\in X:\forall\,x\in G_n\  , \chi(x)=1\}$ be the
annihilator of the group~$G_n$. Each annihilator~$ G_n^\bot$ is
a~group with respect to multiplication, and the subgroups~$
G_n^\bot$ form an~increa\-sing sequence
\begin{equation}
\label{eq1.03} \cdots\subset G_{-n}^\bot\subset\cdots\subset
G_0^\bot \subset G_1^\bot\subset\cdots\subset
G_n^\bot\subset\cdots
\end{equation}
with
$$
\bigcup_{n=-\infty}^{+\infty} G_n^\bot=X \quad {and} \quad
\bigcap_{n=-\infty}^{+\infty} G_n^\bot=\{1\},
$$
the quotient group $ G_{n+1}^\bot/ G_n^\bot$ having order~$p$.
The group of characters~$X$ is a zero-dimensional group with a
basic chain \ref{eq1.03}. The group  may be supplied with the
topology using the chain of subgroups~\ref{eq1.03}, the family of
the cosets $ G_n^\bot\cdot\chi$, $\chi\in X$, being taken
as~a~base of the topology. The collection of such cosets, along
with the empty set, forms the~semiring~$\mathscr{X}$. Given
a~coset $ G_n^\bot\cdot\chi$, we define a~measure~$\nu$ on it by
$\nu( G_n^\bot\cdot\chi)=\nu( G_n^\bot)= p^n$. The measure~$\nu$ can be
extended onto the $\sigma$-algebra of measurable sets in the
standard way. One then forms the absolutely convergent integral
$\displaystyle\int_XF(\chi)\,d\nu(\chi)$ using this measure.

The value~$\chi(g)$ of the character~$\chi$ at an element $g\in G$
will be denoted by~$(\chi,g)$. The Fourier transform~$\widehat f$
of an~$f\in L_2( G)$  is defined~as follows
$$
\widehat f(\chi)=\int_{ G}f(x)\overline{(\chi,x)}\,d\mu(x)=
\lim_{n\to+\infty}\int_{ G_{-n}}f(x)\overline{(\chi,x)}\,d\mu(x),
$$
with the limit being in the norm of $L_2(X)$. For any~$f\in L_2(G)$,
the inversion formula is valid
$$
f(x)=\int_X\widehat f(\chi)(\chi,x)\,d\nu(\chi)
=\lim_{n\to+\infty}\int_{ G_n^\bot}\widehat
f(\chi)(\chi,x)\,d\nu(\chi);
$$
here the limit also signifies the convergence in the norm of~$L_2(
G)$. If $f,g\in L_2( G)$ then the Plancherel formula is valid \cite{AVDR}
$$
\int_{ G}f(x)\overline{g(x)}\,d\mu(x)= \int_X\widehat
f(\chi)\overline{\widehat g(\chi)}\,d\nu(\chi).
$$
\goodbreak

Provided with this topology, the group of characters~$X$ is
a~zero-dimensional locally compact group; there is, however,
a~dual situation: every element $x\in G$ is a~character of the
group~$X$, and~$ G_n$ is the annihilator of the group~$ G_n^\bot$.
We will denote the union of disjoint sets $E_j$  by $\bigsqcup E_j$.

 For any $n\in \mathbb Z$ we choose a character $r_n\in  G_{n+1}^{\bot}\backslash G_n^{\bot}$
 and fixed it. The collection of functions $(r_n)_{n\in \mathbb Z}$ is called a Rademacher system. Any character $\chi$ can be rewritten as a product
 $$\chi=\prod_{j=-m}^{+\infty} r_j^{\alpha_j},\ \alpha_j=\overline{0,p-1}.$$
 Then the Monna mapping $\lambda':X\to[0,+\infty)$ can be defined as
 $$\lambda'(\chi)=\sum_{j=-\infty}^{m}\alpha_{j}p^{j}.$$
 It is evident that $\lambda'(G^\bot_n)=[0,p^n]$.
  Let us denote
  $$
    H_0=\{h\in G: h=a_{-1}g_{-1}\dot+a_{-2}g_{-2}\dot+\dots \dot+ a_{-s}g_{-s}, s\in \mathbb
    N,\ a_j=\overline{0,p-1}\},
  $$
  $$
    H_0^{(s)}=\{h\in G: h=a_{-1}g_{-1}\dot+a_{-2}g_{-2}\dot+\dots \dot+
    a_{-s}g_{-s},\ a_j=\overline{0,p-1}
    \},s\in \mathbb N.
  $$
  Under the Monna mapping $\lambda(H_0)= \mathbb N_0=\mathbb N\bigsqcup \{0\}$ and $\lambda(H_0^{(s)})=\mathbb N_0\bigcap [0,p^{s-1}].$ Thus the set $H_0$ is an analog of the set $\mathbb N_0$.
  \begin{definition}
 We define the mapping ${\cal A}\colon G\to G$ by
 ${\cal A}x:=\sum_{n=-\infty}^{+\infty}a_ng_{n-1}$, where
 $x=\sum_{n=-\infty}^{+\infty}a_ng_n\in G$. As any element $x\in G$ can
 be uniquely expanded as~$x=\sum a_ng_n$, the mapping ${\cal A}\colon G\to
 G$ is one-to-one onto. The mapping~${\cal A}$ is called
 a  dilation operator if~${\cal A}(x\dot+ y)={\cal A}x\dot + {\cal A}y$ for all
 $x,y\in G$.
 \end{definition}
 We note that if $ G$~is a~Vilenkin group ($p\cdot g_n=0$)
 or is the group of all $p$-adic numbers ($p\cdot g_n=g_{n+1}$),
 then~${\cal A}$ is an~additive operator and hence a dilation operator.
   Moreover, if there exist fixed numbers $c_1,c_2,\dots,c_{\tau} =\overline{0,p-1}$ such that \begin{equation}
\label{eq1.04}
  pg_n=c_1g_{n+1}\dot+ c_2g_{n+2}\dot+
 \cdots\dot+ c_{\tau}g_{n+\tau},
 \end{equation}
  then the operator  ${\cal A}$ will be additive.
  We will assume that the condition (\ref{eq1.04}) is satisfied.
     By definition, put $(\chi {\cal A},x)=(\chi, {\cal
 A}x)$. It is also clear that
 ${\cal A} g_n= g_{n-1}, r_n{\cal A} = r_{n+1}$,\
 ${\cal A} G_n= G_{n-1}, G_n^\bot{\cal A}= G_{n+1}^\bot$.

  \begin{lemma}[\cite{LS2}]
 For any zero-dimensional group\\
 1) $\int\limits_{G_0^\bot}(\chi,x)\,d\nu(\chi)={\bf 1}_{G_0}(x)$,
 2) $\int\limits_{G_0}(\chi,x)\,d\mu(x)={\bf 1}_{G_0^\bot}(\chi)$.\\
   3) $\int\limits_{G_n^\bot}(\chi,x)\,d\nu(\chi)=p^n{\bf
  1}_{G_n}(x)$,
  4) $\int\limits_{G_n}(\chi,x)\,d\mu(x)=\frac{1}{p^n}{\bf
  1}_{G_n^\bot}(\chi)$.
  \end{lemma}
\begin{lemma}[\cite{LS2}]
Let  $\chi_{n,s}=r_n^{\alpha_n}r_{n+1}^{\alpha_{n+1}}\dots
r_{n+s}^{\alpha_{n+s}}$ be a character which does not belong to
$G_n^\bot$. Then
$$
\int\limits_{G_n^\bot\chi_{n,s}}(\chi,x)\,d\nu(\chi)=p^n(\chi_{n,s},x){\bf
1}_{G_n}(x).
$$
\end{lemma}
\begin{lemma}[\cite{LS2}]
Let
$h_{n,s}=a_{n-1}g_{n-1}\dot+a_{n-2}g_{n-2}\dot+\dots\dot+a_{n-s}g_{n-s}\notin
G_n$. Then
$$
\int\limits_{G_n\dot+h_{n,s}}(\chi,x)\,d\mu(x)=\frac{1}{p^n}(\chi,h_{n,s}){\bf
1}_{G_n^\bot}(\chi).
$$
\end{lemma}
 \begin{definition}[\cite{LS2}]
Let $M,N\in\mathbb N$.
We denote by  ${\mathfrak D}_M(G_{-N})$ the set of functions
 $f\in L_2(G)$ such that 1) ${\rm supp}\,f\subset G_{-N}$, and 2)
 $f$ is constant on cosets $G_M\dot+g$. The class ${\mathfrak
 D}_{-N}(G_{M}^\bot)$ is defined similarly.
 \end{definition}

 \begin{lemma}\label{lm.1.1}
 For any fixed  $\alpha_0,\alpha_1,...,\alpha_s =0,1,...,p-1$  the set  $H_0$ is an orthonormal basis in
   $L_2(G_0^\bot r_0^{\alpha_0}r_1^{\alpha_1}...r_s^{\alpha_s})$.
 \end{lemma}
 {\it Proof.} It is known
   (\cite{LSF})that the set  $H_0$ is  an orthonormal system in  $L_2(G_0^\bot)$.
 Moreover,  $H_0$ is an orthonormal basis in $L_2(G_0^\bot)$. Indeed, for any  $\nu\in \mathbb N$ functions
 $
 h_j=a_{-1}g_{-1}\dot+a_{-2}g_{-2}\dot+...\dot+a_{-\nu}g_{-\nu}
 $
 are constant on cosets
 \begin{equation}\label{eq1.05}
 G_{-\nu}^\bot r_{-\nu}^{\alpha_{-\nu}}r_{-\nu+1}^{\alpha_{-\nu+1}}...r_{-1}^{\alpha_{-1}}
 \end{equation}
   and orthogonal on the set    $G_0^\bot$.
    Therefore, any step function that is constant on cosets (\ref{eq1.05})
   can be uniquely represented us linear combination of elements
    $$
 h_j=a_{-1}g_{-1}\dot+a_{-2}g_{-2}\dot+...\dot+a_{-\nu}g_{-\nu}.
 $$
  So the system  $H_0$ is an orthonormal basis in $L_2(G_0^\bot)$.
 Let us show that $H_0$ is an orthonormal basis in  $L_2(G_0^\bot r_0^{\alpha_0}r_1^{\alpha_1}...r_s^{\alpha_s})$. For any $h_1,h_2\in H_0$
 $$
   \int_{G_0^\bot r_0^{\alpha_0}r_1^{\alpha_1}...r_s^{\alpha_s}}(\chi,h_1)\overline{(\chi,h_2)}d\nu(\chi)=
   \int_X{\bf 1}_{G_0^\bot r_0^{\alpha_0}r_1^{\alpha_1}...r_s^{\alpha_s}}(\chi)(\chi,h_1)\overline{(\chi,h_2)}d\nu(\chi)=
  $$
  $$
   =
   \int_X{\bf 1}_{G_0^\bot r_0^{\alpha_0}r_1^{\alpha_1}...r_s^{\alpha_s}}(\chi r_0^{\alpha_0}r_1^{\alpha_1}...r_s^{\alpha_s})(\chi r_0^{\alpha_0}r_1^{\alpha_1}...r_s^{\alpha_s},h_1)\overline{(\chi r_0^{\alpha_0}r_1^{\alpha_1}...r_s^{\alpha_s},h_2)}d\nu(\chi)=
  $$
  $$ =
   \int_X{\bf 1}_{G_0^\bot }(\chi)(\chi r_0^{\alpha_0}r_1^{\alpha_1}...r_s^{\alpha_s},h_1)\overline{(\chi r_0^{\alpha_0}r_1^{\alpha_1}...r_s^{\alpha_s},h_2)}d\nu(\chi)=
  $$
  $$ =
   (r_0^{\alpha_0}r_1^{\alpha_1}...r_s^{\alpha_s},h_1)
   \overline{(r_0^{\alpha_0}r_1^{\alpha_1}...r_s^{\alpha_s},h_2)}
   \int_X{\bf 1}_{G_0^\bot }(\chi)(\chi ,h_1)\overline{(\chi,h_2)}d\nu(\chi)=
  $$
  $$ =
   (r_0^{\alpha_0}r_1^{\alpha_1}...r_s^{\alpha_s},h_1)
   \overline{(r_0^{\alpha_0}r_1^{\alpha_1}...r_s^{\alpha_s},h_2)}
   \int_{G_0^\bot} (\chi ,h_1)\overline{(\chi,h_2)}d\nu(\chi)=\delta_{h_1,h_2}.
  $$
 This means that $H_0$ is an orthonormal system.   Let us show that $H_0$ is a basis in  $L_2(G_0^\bot r_0^{\alpha_0}r_1^{\alpha_1}...r_s^{\alpha_s})$.
    For any $\nu\in \mathbb N$ functions
 $
 h_j=a_{-1}g_{-1}\dot+a_{-2}g_{-2}\dot+...\dot+a_{-\nu}g_{-\nu}
 $
 are constant on cosets
  \begin{equation}\label{eq1.07}
 G_{-\nu}^\bot r_{-\nu}^{\alpha_{-\nu}}r_{-\nu+1}^{\alpha_{-\nu+1}}...r_{-1}^{\alpha_{-1}}
 r_0^{\alpha_0}r_1^{\alpha_1}...r_s^{\alpha_s}
 \end{equation}
  and orthogonal on the set   $G_0^\bot  r_0^{\alpha_0}r_1^{\alpha_1}...r_s^{\alpha_s}$.
  Therefore, any step function that is constant on cosets (\ref{eq1.07})
   can be uniquely represented us linear combination of elements
    $$
 h_j=a_{-1}g_{-1}\dot+a_{-2}g_{-2}\dot+...\dot+a_{-\nu}g_{-\nu}.
 $$
 This means that $H_0$ is an orthonormal basis.
    $\square$

  \begin{lemma}\label{lm.1.2} Let $s\in \mathbb N$.
 For any fixed $\alpha_{-1},...,\alpha_{-s}=\overline{0,p-1}$ the family   $p^{\frac{s}{2}}\mathcal{A}^sH_0$
  is  an orthonormal basis in  $L_2(G_{-s}^\bot r_{-s}^{\alpha_{-s}}...r_{-1}^{\alpha_{-1}} )$.
 \end{lemma}
   {\it Proof.} Suppose $\tilde{h}_{1},\tilde{h}_{2}\in \mathcal{A}^s H_0$;  then
  $$
  \int_{G_{-s}^\bot r_{-s}^{\alpha_{-s}}...r_{-1}^{\alpha_{-1}}}(\chi,\tilde{h}_{1})\overline{(\chi,\tilde{h}_{2})}d\nu(\chi)=
  \int_{G_{-s}^\bot r_{-s}^{\alpha_{-s}}...r_{-1}^{\alpha_{-1}}}(\chi\mathcal{A}^s,\mathcal{A}^{-s}\tilde{h}_{1})
  \overline{(\chi\mathcal{A}^s,\mathcal{A}^{-s}\tilde{h}_{2})}d\nu(\chi)=
  $$
  $$
  \int_X{\bf 1}_{G_{-s}^\bot r_{-s}^{\alpha_{-s}}...r_{-1}^{\alpha_{-1}}}(\chi)(\chi\mathcal{A}^s,\mathcal{A}^{-s}\tilde{h}_{1})
  \overline{(\chi\mathcal{A}^s,\mathcal{A}^{-s}\tilde{h}_{2})}d\nu(\chi)=
   $$
  $$
  =\frac{1}{p^s}\int_X{\bf 1}_{G_{-s}^\bot r_{-s}^{\alpha_{-s}}...r_{-1}^{\alpha_{-1}}}(\chi\mathcal{A}^{-s})(\chi,\mathcal{A}^{-s}\tilde{h}_{1})
  \overline{(\chi,\mathcal{A}^{-s}\tilde{h}_{2})}d\nu(\chi)=
   $$
   $$
  =\frac{1}{p^s}\int_{G_{0}^\bot r_0^{\alpha_{-s}}r_1^{\alpha_{-s+1}}... r_{s-1}^{\alpha_{-1}}   }(\chi,\mathcal{A}^{-s}\tilde{h}_{1})
  \overline{(\chi,\mathcal{A}^{-s}\tilde{h}_{2})}d\nu(\chi)=\frac{1}{p^s}\delta_{\tilde{h}_1,\tilde{h}_2}.
   $$
   So, the system  $p^{\frac{s}{2}}\mathcal{A}^sH_0$ is orthonormal.   A basis property  is proved as before. $\square$


 For a given function $\varphi\in  {\mathfrak D}_M(G_{-N})$ , we define  subspaces $V_n\subset L_2(G)$ generated by $\varphi$ as
  $$
 V_n=\overline{{\rm span}\{\varphi(\mathcal{A}^n\cdot \dot-h),h\in H_0\}}, \ n\in\mathbb Z.
 $$

 We say that sequence of subspaces $\{V_n\}$ forms a {\it multiresolution analysis} (MRA) for $L_2(G)$, if the folloving conditions are satisfied
 \begin{equation}\label{eq1.08}
 V_n\subset V_{n+1}, n\in \mathbb Z,
 \end{equation}
 \begin{equation}\label{eq1.09}
 \overline{\cup_n V_n}=L_2(G),
 \end{equation}
  \begin{equation}\label{eq1.1}
 \cap_n V_n=\{0\}
 \end{equation}

A function $\varphi \in L_2(G)$  is called refinable if
  \begin{equation}                                      \label{eq1.2}
   \varphi(x)=p\sum_{h\in H_0}\beta_h\varphi({\cal
   A}x\dot-h),
   \end{equation}
      for some sequence  $(\beta_h)\in \ell^2$. The equality   (\ref{eq1.2}) is called the refinement  equation. In Fourier domain, the equality   (\ref{eq1.2})  can be rewritten as
  $$
 \hat\varphi(\chi)=m_0(\chi)\hat\varphi(\chi{\cal A}^{-1}),
 $$
 where
 \begin{equation}                                           \label{eq1.3}
 m_0(\chi)=\sum_{h\in
 H_0}\beta_h\overline{(\chi{\cal A}^{-1},h)}
 \end{equation}
   is a mask of  (\ref{eq1.2}).
   \begin{lemma}
    If a refinable function  $\varphi\in \mathfrak{D}_{G_{M}}(G_{-N}),\ M,N\in \mathbb N$, then its refinement equation is
      \begin{equation}                                      \label{eq1.4}
   \varphi(x)=p\sum_{h\in H_0^{(N+1)}}\beta_h\varphi({\cal
   A}x\dot-h),
   \end{equation}
   \end{lemma}
   {\bf Proof.} Let us write  $\varphi$  in the form
   \begin{equation}                                      \label{eq1.5}
   \varphi(x)=p\sum_{h\in H_0^{(N+1)}}\beta_h\varphi({\cal
   A}x\dot-h)+p\sum_{h\notin H_0^{(N+1)}}\beta_h\varphi({\cal
   A}x\dot-h).
   \end{equation}
   1)Let $x\in G_{-N}$. Then   ${\cal A} x\in G_{-N-1}$. If $h\notin H_0^{(N+1)}$,
   then
   $$
   h=a_{-1}g_{-1}\dot+...\dot+ a_{-N-1}g_{-N-1}\dot+  a_{-N-2}g_{-N-2}\dot+ ...,
   $$
    where $a_{-N-2}g_{-N-2}\dot+ ...\neq 0$. It  follows  ${\cal A}x \dot-h \notin G_{-N-1}$ and
       $$
   \sum_{h\notin H_0^{(n+1)}}\beta_h\varphi({\cal A}x\dot-h)=0.
   $$
   2)Let $x\notin G_{-N}$. Then $\varphi(x)=0$ and  ${\cal A} x\notin G_{-N-1}$.
   Therefore
   $$
   {\cal A}x =\dots \dot+a_{-1}g_{-1}\dot+...\dot+ a_{-N-1}g_{-N-1}\dot+  a_{-N-2}g_{-N-2}\dot+ ...,
   $$
    where $a_{-N-2}g_{-N-2}\dot+ ...\neq 0$.    If $h\in H_0^{(N+1)}$, then
      $$
   {\cal A}x\dot-h =\dots \dot+b_{-1}g_{-1}\dot+...\dot+ b_{-N-1}g_{-N-1}\dot+  a_{-N-2}g_{-N-2}\dot+ ...,
   $$
       where $a_{-N-2}g_{-N-2}\dot+ ...\neq 0$. It follows  $\varphi({\cal A}x\dot-h)=0$ and $\sum_{h\notin H_0^{(n+1)}}\beta_h\varphi({\cal A}x\dot-h)=0. \quad \square$


 If the shift system
 $(\varphi(x\dot-h))_{h\in H_0}$ form an  orthonormal basis in $V_0$, then MRA $(V_n)$ is called orthogonal.  Orthogonal MRA is used to construct orthogonal affine systems that form a basis of $L_2( G)$.
  \begin{theorem}[ \cite{LSF}.] Let $\varphi\in \mathfrak{D}_{G_{0}}(G_{-N})$ be a   refinable function and
  $|\hat\varphi(\chi)|={\bf 1}_{G_0^\bot}(\chi)$. Then $\varphi$ generates an orthogonal MRA .
 \end{theorem}
 We can choose a mask $m_0(\chi)\in \mathfrak{D}_{G_{-N}^\bot}(G_{0}^\bot)$ so that   $ m_0(G_{-N}^\bot)=1, |m_0(\chi)|={\bf 1}_{G_0^\bot}(\chi) $.  Then the corresponding scaling function $\varphi$ generates an orthogonal MRA. In this case, orthogonal wavelets $\hat\psi_\ell(\chi)$ are defined by the equalities
 \begin{equation}                                      \label{eq1.6}
   \psi^{(\ell)}(x)=\sum_{h\in H_0}\beta_h^{(\ell)}\varphi({\cal
   A}x\dot-h).
   \end{equation}
 In the Fourier domain, (\ref{eq1.6}) can be written  as
 \begin{equation}\label{eq1.7}
 \hat\psi_\ell(\chi)=\hat\varphi(\chi A^{-1})m_\ell(\chi), \quad (\ell=1,2,\dots,p-1),
  \end{equation}
   where
    $m_\ell(\chi)=m_0(\chi r_0^{-\ell})\;\;(\ell=1,2,\dots,p-1)$.
 The system \\$(\psi_\ell(x\dot- h))_{h\in H_0, \ell=1,2,\dots,p-1}$ is an orthogonal basis of the orthogonal complement $V_1\ominus V_0=\{x\in V_1:x\bot V_0\}$  \cite{LSF}.

If the shifts $(\varphi(x\dot-h))_{h\in H_0}$  are not orthogonal, then one can try to choose the functions $\psi^{(\ell)}(x)$  so that for any $ f\in L_2(G)$
   $$
   f(x)=\sum_{\ell =1}^r\sum_{n\in \mathbb Z}\sum_{h\in H_0}(f,\psi^{(\ell)} ({\cal A}^n \cdot\dot - h))\psi^{(\ell)} ({\cal A}^n x\dot - h).
   $$
    Such  system is called Parseval   wavelet frame or tight wavelet frame.

\section{Refinable functions in  p-adic groups}
In this section we will construct a  refinable function in the group $\mathfrak{G}$ with the condition $pg_n=g_{n+1}$.  It means that $\mathfrak{G}$ is addidive group of the field $\mathbb Q_p $ of all $p$-adic numbers. In this case $(r_n,g_m)=e^{\frac{2\pi i}{p^{n-m+1}}}$.

Let $M,N\in \mathbb N$.  We want to construct a refinable function
 $\varphi\in \mathfrak D_{\mathfrak{G}_M}(\mathfrak{G}_{-N})$, i.e.
$\hat{\varphi}\in \mathfrak D_{\mathfrak{G}_{-N}^\bot}(\mathfrak{G}_{M}^\bot)$.
Let us write the refinable equation in the Fourier domain

 \begin{equation}\label{eq2.1}
 \hat\varphi_(\chi)=\hat\varphi(\chi A^{-1})m_0(\chi).
  \end{equation}
 The mask
  \begin{equation} \label{eq2.2}
  m_0(\chi)=\sum_{h\in
  H_0^{(N+1)}}\beta_h\overline{(\chi,A^{-1}h)}
 \end{equation}
 is constant on cosets  $\mathfrak{G}_{-N}^\bot r_{-N}^{\alpha_{-N}}...r_{-N+s}^{\alpha_{-N+s}}$ .
 Let $m_0(\mathfrak{G}_{-N}^\bot )=1$. Then
$$
\hat{\varphi}(\chi)=m_0(\chi)m_0(\chi\mathcal A^{-1})...m_0(\chi\mathcal A^{-N-M}).
$$
Denote by
 \begin{equation} \label{eq2.3}
 m_0(\mathfrak{G}_{-N}^\bot r_{-N}^{\alpha_{-N}}r_{-N+1}^{\alpha_{-N+1}}...r_{0}^{\alpha_{0}}...r_{M}^{\alpha_{M}})=: \lambda_{ \alpha_{-N}\alpha_{-N+1}...\alpha_{0}...\alpha_{M}}=\lambda_m,
 \end{equation}
 where
 $$
 m= \alpha_{-N}+\alpha_{-N+1}p+...+\alpha_{0}p^N +...+\alpha_M p^{N+M} ,
 $$
  the values of the mask on  $\mathfrak{G}_{M+1}^\bot$. Since $(\chi \mathcal{A}^{-1},h)$ is constant on cosets
 $$
 \mathfrak{G}_{-N}^\bot r_{-N}^{\alpha_{-N}}r_{-N+1}^{\alpha_{-N+1}}...r_{M}^{\alpha_{M}}
$$
we have
$$
(\chi\mathcal{A}^{-1}, h)=(\mathfrak{G}_{-N}^\bot r_{-N}^{\alpha_{-N}}...r_{M}^{\alpha_{M}},
a_{-1}g_{0}\dot+a_{-2}g_{-1}\dot+...+a_{-N-1}g_{-N})=,
$$
$$
=( r_{-N}^{\alpha_{-N}}...r_{M}^{\alpha_{M}},
a_{-1}g_{0}\dot+a_{-2}g_{-1}\dot+...+a_{-N-1}g_{-N})=,
$$
$$
=\prod\limits_{\nu=0}^{-N} \prod\limits_{k=-N}^M ( r_{k},g_{\nu})^{\alpha_ka_{\nu-1}}=
\prod\limits_{\nu=0}^{-N} \prod\limits_{k=-N}^M e^{\frac{2\pi i}{p^{k-\nu+1}}\alpha_k a_{\nu}}.
$$

Therefore
$$
m_0(\chi)=\sum\limits_{a_{-1},...,a_{-N-1}}\beta_{a_{-1},...,a_{-N-1}}
e^{-\frac{2\pi i}{p}\sum\limits_{\nu=0}^{-N}a_{\nu-1}p^\nu\sum\limits_{k=-N}^M \frac{\alpha_k}{p^k} }=
$$
$$
=\sum\limits_{a_{-1},...,a_{-N-1}}\beta_{a_{-1},...,a_{-N-1}}
e^{-\frac{2\pi i}{p}\sum\limits_{\nu=0}^{-N}a_{\nu-1}p^{\nu+N}\sum\limits_{k=-N}^M \alpha_k p^{-k
-N} }.
$$
Denote

$$
q_m=e^{-\frac{2\pi i}{p}}\sum\limits_{k=-N}^{M}\alpha_k p^{-N-k}, \ n=a_{-N-1}+a_{-N}p+...+a_{-1}p^N.
$$
and write the equality  (\ref{eq2.2})  in the form

\begin{equation}\label{eq2.4}
 \left(\begin{array}{rrrr}
q_{0}^0&q_{0}^1&...&q_{0}^{p^{N+1}-1}\\
q_{1}^0&q_{1}^1&...&q_{1}^{p^{N+1}-1}\\
...&...&...&...\\
q_{p^{N+1}-1}^0&q_{p^{N+1}-1}^1&...&q_{p^{N+1}-1}^{p^{N+1}-1}\\
...&...&...&...\\
q_{p^{M+N+1}-1}^0&q_{p^{M+N+1}-1}^1&...&q_{p^{M+N+1}-1}^{p^{N+1}-1}\\

\end{array}
\right)
\left(\begin{array}{c}
\beta_0\\
\beta_1\\
\beta_2\\
\vdots\\
\beta_{p^{N+1}-1}
\end{array}
\right)=
\left(\begin{array}{c}
\lambda_0\\
\lambda_1\\
\lambda_2\\
\vdots\\
\lambda_{p^{M+N+1}-1}\\
\end{array}
\right)
 \end{equation}
We need to find $\lambda_m$ and $\beta_n$ so that
$$
\hat{\varphi}(\chi)=m_0(\chi)m_0(\chi\mathcal A^{-1})...m_0(\chi\mathcal A^{-N-M})=0
$$
on the set $\mathfrak{G}_{M+1}^\bot\setminus \mathfrak{G}_{M}^\bot$.
To find $\lambda_m$  and $\beta_n$ we construct a rooted $p$-{\it adic} mask tree $T=T(m_0)$ in the following way.
For any $m\in \mathbb N: p^{M+N}\le m\le p^{M+N+1}-1$ we construct the path
$$
\lambda_m\rightarrow \,\lambda_{m \,{\rm div}\, p} \rightarrow \lambda_{m \,{\rm div}\, p^2}\dots \rightarrow \lambda_{m \,{\rm div}\, p^{M+N}}\rightarrow \lambda_0=1.
$$
 from the leaf $\lambda_m$ to the root $\lambda_0$ of the  tree $T(m_0)$. It is clear that $H=M+N$ is a high of this $p$-{\it adic} tree $T$. (See Figure 2 for graph of $T$). In this tree, the numbers  $\lambda_{j_s}: p^{s-1}\le j_s\le p^{s}-1, s\ge 0$ form the s-th level. The set of all products $\lambda_{m}\lambda_{m \,{\rm div} \, p}...\lambda_0$ coincides with the set of all values of the function $\hat{\varphi}(\chi)$ on the set $\mathfrak{G}_{M+1}^\bot\setminus \mathfrak{G}_{M}^\bot$. On each path
$$
\lambda_m\rightarrow \,\lambda_{m \,{\rm div}\, p} \rightarrow \lambda_{m \,{\rm div}\, p^2}\dots \rightarrow \lambda_{m \,{\rm div}\, p^{M+N}}
$$
 we select one node and place zero there. Denote $\Lambda_0 (T))=\{\lambda_\nu =0\}$.

\unitlength=0.80mm
  \begin{picture}(240,80)
 \small

  \put(150,48){\line (-1,1){16}}
     \put(151,48){\line (0,1){16}}
     \put(152,48){\line (1,1){16}}
     \put(125,66){$\lambda_{p^{N+M+1}-p}$}
      \put(160,66){$\lambda_{p^{N+M+1}-1}$}

     \put(58,48){\line (-1,1){16}}
     \put(59,48){\line (0,1){16}}
     \put(60,48){\line (1,1){16}}
     \put(38,66){$\lambda_{p^{N+M}}$}
      \put(70,66){$\lambda_{p^{N+M}+p-1}$}
      \multiput(51,64)(3,0){6}{$\cdot$}

    \multiput(60,46)(3,0){31}{$\cdot$}
   \put(130,25){\line (-1,1){16}}
     \put(131,25){\line (0,1){16}}
     \put(132,25){\line (1,1){16}}
     \put(110,43){$\lambda_{p^{N}-p}$}
      \put(146,43){$\lambda_{p^{N}-1}$}
      \multiput(124,41)(3,0){6}{$\cdot$}
      \multiput(74,41)(3,0){4}{$\cdot$}

  \put(124,19){$\lambda_{p-1}$}
  \put(81,19){$\lambda_{1}$}
   \put(105,19){$\lambda_{\nu}$}
   \put(80,25){\line (-1,1){16}}
     \put(81,25){\line (0,1){16}}
     \put(82,25){\line (1,1){16}}
     \put(60,43){$\lambda_{p^{N-1}}$}
      \put(86,43){$\lambda_{p^{N-1}+p-1}$}

  \multiput(81,22)(3,0){17}{$\cdot$}

 \put(100,-2){$\lambda_{0}=1$}
  \put(100,2){\line (-1,1){16}}
  \put(106,2){\line (0,1){16}}
  \put(110,2){\line (1,1){16}}
  \put(-5,66){$\bf \mathfrak{G}_{M+1}^\bot \setminus \mathfrak{G}_{M}^\bot $:}
  \put(-5,43){$\bf \mathfrak{G}_{0}^\bot \setminus \mathfrak{G}_{-1}^\bot $:}
  \put(-5,19){$\bf \mathfrak{G}_{-N+1}^\bot \setminus \mathfrak{G}_{-N}^\bot $:}
  \put(-5,-2){$\bf \mathfrak{G}_{-N}^\bot  $:}

    \end{picture}\\

\hskip4cm Figure 2.  The graph of the tree $T=T(m_0)$.

\begin{theorem}\label{th2.01}
 1) Let $\sharp\Lambda_0(T)) = p^{N+1}-1$. Then  the corresponding  values $\lambda_\nu\in \Lambda_0(T)$ determine the mask of the refinable  function. If $\lambda_\nu\notin \Lambda_0(T) $  then $\lambda_\nu\neq 0$.\\
 2) Let $\sharp\Lambda_0(T) < p^{N+1}-1$. Then  the corresponding  values $\lambda_\nu\in \Lambda_0(T)$ determine the mask of the refinable  function.\\
 3) If $\sharp\Lambda_0(T)\ge p^{N+1}$ then  the corresponding  values $\lambda_\nu\in \Lambda_0(T)$ do not define the mask of the refinable  function.
 \end{theorem}
{\bf Proof.}1)First we consider the case $\sharp\Lambda_0(T)= p^{N+1}-1$.
 Let $\Lambda_0(T)=(\lambda_{\nu_j})_{j=1}^{ p^{N+1}-1}.$ Let us choose equations with numbers $0,\nu_1,\nu_2,\dots ,\nu_{p^{N+1}-1}$ from system (\ref{eq2.4}). Then we get a system with Vandermonde determinant. Solving this system, we find the coefficients $\beta_0,\beta_1,...,\beta_{p^{N+1}-1}$. After that, we find the remaining values $\lambda_j\notin \Lambda_0(T) $.  Let us show that these values $\lambda_j\neq 0 $. Suppose that this is false.  Then there exists $\lambda_{\nu_0}=0, \ \lambda_{\nu_0}\in \Lambda_0(T)) $.  Let us choose equations with numbers $\nu_0,\nu_1,\nu_2,\dots ,\nu_{p^{N+1}-1}$ from system (\ref{eq2.4}). Then we get a system with Vandermonde determinant. Solving this system, we find the coefficients $\beta_0=\beta_1=...=\beta_{p^{N+1}-1}=0$ . It follows that $\lambda_0=0$, that is impossible.\\
 2)Let $\sharp\Lambda_0(T)) < p^{N+1}-1$. Add  values $\lambda_\nu =0 $ to the set $\Lambda_0(T)$ and use the first point.\\
 3)Let $\sharp\Lambda_0(T) > p^{N+1}-1$. Let $\Lambda_0(T)=(\lambda_{\nu_j})_{j=1}^{ p^{N+1}+l},\, l>-1.$ Let us choose equations with numbers $\nu_1,\nu_2,\dots ,\nu_{p^{N+1}}$ from system (\ref{eq2.4}). Then we get a system with Vandermonde determinant. Solving this system, we find the coefficients $\beta_0=\beta_1=...=\beta_{p^{N+1}-1}=0$. It follows that $\lambda_0=0$, that is impossible. $\square$\\

 So we have some method  to construct step scaling functions
 $\varphi\in \mathfrak D_{\mathfrak{G}_{M}}(\mathfrak{G}_{-N})$.\\

\section{Constructing  tight wavelet frames without the principle of unitary extension   }
For the function $\varphi\in L_2(\mathfrak{G})$ we will use the standard notation
$$
\varphi_{n,h}=p^{\frac{n}{2}}\varphi(\mathcal{A}^n\cdot \dot-h), \quad h\in H_0, n\in \mathbb Z.
$$
Let the mask  $m_0({\chi})$ built on a tree $T$ and $\varphi\in \mathfrak D_{\mathfrak{G}_{M}}(\mathfrak{G}_{-N})$ be the corresponding refinable  function.
  We want to find masks $m_j, j=\overline{1,q}$ and corresponding functions $\psi^{(j)}$ that generate a tight wavelet frame. We will use lemma  \ref{Lm3.1} and theorem \ref{Th3.1} which are valid for any zero-dimensional group $G$.
\subsection{Sufficient conditions in terms of masks. }
\begin{lemma}\label{Lm3.1}
Let the mask $m_j \ (j=1,...,q)$ satisfy the following conditions:
1)$
\hat{\psi}^{(j)}(\chi)=\hat{\varphi}(\chi\mathcal{A}^{-1})m_j(\chi)=1
 $   on the coset
 $   G^\bot_{-s}r_{-s}^{\gamma_{-s}}
... r_{-1}^{\gamma_{-1}}r_{0}^{\gamma_{0}}...r_{u}^{\gamma_{u}}, \quad s=s(j), s\le N
$
 and\\
 2) $m_j(\chi)=0$ outside the coset $   G^\bot_{-s}r_{-s}^{\gamma_{-s}}
... r_{-1}^{\gamma_{-1}}r_{0}^{\gamma_{0}}...r_{u}^{\gamma_{u}}$.\\
Then
$$
\sum_{h\in H_0} |c_{n,h}^{(j)}(f)|^2=
\sum_{h\in H_0} |(\psi^{(j)}_{n,h},f)|^2=
\int_{G_{n-s}^\bot r_{n-s}^{\gamma_{-s}}...r_{n-1}^{\gamma_{-1}}r_{n-0}^{\gamma_{0}}...r_{n+u}^{\gamma_{u}}}|\hat{f}(\chi) |^2d\nu(\chi)
$$
\end{lemma}
{\bf Proof.} By the properties of the Fourier transform, we have
\begin{equation}\label{Eq3.1}
c_{n,h}^{(j)}(f)=
p^{\frac{n}{2}}\int_{G}f(x)\overline{\psi^{(j)}(\mathcal{A}^nx\dot-h)}d\mu(x)=
p^{\frac{n}{2}}\int_{X}\hat{f}(\chi)\overline{\hat{\psi}^{(j)}_{\mathcal{A}^n\cdot \dot-h}(\chi)}d\nu(\chi).
\end{equation}
Calculate the Fourier transform
    $$
    \hat{\psi}^{(j)}_{\mathcal{A}^n\cdot \dot-h}(\chi)=\int_{G}\psi^{(j)}(\mathcal{A}^nx \dot -h)\overline{(\chi,x)}d\mu(x)=
    $$
    $$
    \frac{1}{p^n}\int_{G}\psi^{(j)}(x \dot -h)\overline{(\chi\mathcal{A}^{-n},x)}d\mu(x)=
    $$
     $$
    =\frac{1}{p^n}\int_{G}\psi^{(j)}(x) \overline{(\chi\mathcal{A}^{-n},x\dot+ h)}d\mu(x)=
    \frac{1}{p^n}\hat{\psi}^{(j)}(\chi\mathcal{A}^{-n})\overline{(\chi\mathcal{A}^{-n},h)}.
    $$
    Substituting    in  (\ref{Eq3.1}) , we get
    $$
    c_{n,h}^{(j)}(f)=p^{-\frac{n}{2}}\int_{X}\hat{f}(\chi)
    \overline{\hat{\psi}^{(j)}(\chi\mathcal{A}^{-n})}
    (\chi\mathcal{A}^{-n},h)d\nu(\chi)=
    $$
    $$
    =p^{-\frac{n}{2}}\int_{X}\hat{f}(\chi)\overline{\hat{\varphi}(\chi\mathcal{A}^{-n-1})
    m_j(\chi\mathcal{A}^{-n})}(\chi\mathcal{A}^{-n},h)d\nu(\chi)=
    $$
    $$
    =p^{\frac{n}{2}}\int_{X}\hat{f}(\chi \mathcal{A}^{n})\overline{\hat{\varphi}(\chi\mathcal{A}^{-1})
    m_j(\chi)}(\chi,h)d\nu(\chi)=
    $$
    $$
    p^{\frac{n}{2}}\int_{G_{-s}^\bot r_{-s}^{\gamma_{-s}}...r_{0}^{\gamma_{0}}...r_{u}^{\gamma_{u}}  }\hat{f}(\chi \mathcal{A}^{n})(\chi,h)d\nu(\chi).
    $$
    Let us denote $\tilde H_0^{(s)}=\{\tilde h=a_{-s-1}g_{-s-1}\dot+...\}$,
 $h=a_{-1}g_{-1}\dot+...\dot+a_{-s}g_{-s}\dot+ \tilde h$.
  Sinse
  $(\chi, a_{-1}g_{-1}\dot+...\dot+a_{-s}g_{-s})$  is constant on any coset  $G_{-s}^\bot r_{-s}^{\gamma_{-s}}...r_{0}^{\gamma_{0}}...r_{u}^{\gamma_{u}}$ it follows by lemma  \ref{lm.1.2} that
$$
    \sum_{h\in H_0}|c_{n,h}^{(l)}(f)|^2=
$$
$$
=\sum_{a_{-1}=0}^{p-1}...\sum_{a_{-s}=0}^{p-1}\sum_{\tilde h\in \tilde H_0^{(s)}}
    |p^{\frac{n}{2}}\int_{G_{-s}^\bot r_{-s}^{\gamma_{-s}}...r_{0}^{\gamma_{0}}...r_{u}^{\gamma_{u}}  }\hat{f}(\chi \mathcal{A}^{n})(\chi,\tilde h)d\nu(\chi)|^2=
$$
$$
=p^{s+n}\sum_{\tilde h\in \tilde H_0^{(s)}}
    |\int_{G_{-s}^\bot r_{-s}^{\gamma_{-s}}...r_{0}^{\gamma_{0}}...r_{u}^{\gamma_{u}}  }\hat{f}(\chi \mathcal{A}^{n})(\chi,\tilde h)d\nu(\chi)|^2=
$$
$$
=p^n\int_{G_{-s}^\bot r_{-s}^{\gamma_{-s}}...r_{0}^{\gamma_{0}}...r_{u}^{\gamma_{u}}  }|\hat{f}(\chi \mathcal{A}^{n})|^2d\nu(\chi)=
$$
$$
=p^n\int_X {\bf 1}_{G_{-s}^\bot r_{-s}^{\gamma_{-s}}...r_{0}^{\gamma_{0}}...r_{u}^{\gamma_{u}}}(\chi)|\hat{f}(\chi \mathcal{A}^{n})|^2d\nu(\chi)=
$$
$$
=\int_X {\bf 1}_{G_{-s}^\bot r_{-s}^{\gamma_{-s}}...r_{0}^{\gamma_{0}}...r_{u}^{\gamma_{u}}}(\chi\mathcal{A}^{-n})|\hat{f}(\chi) |^2d\nu(\chi)=
$$
$$
=\int_{G_{-s}^\bot r_{-s}^{\gamma_{-s}}...r_{0}^{\gamma_{0}}..r_{u}^{\gamma_{u}}\mathcal{A}^{n}}|\hat{f}(\chi) |^2d\nu(\chi).
$$
The lemma is proved.  $\square$

\begin{theorem}\label{Th3.1}
  Let $\varphi \in \mathfrak{D}_{G_{M}}(G_{-N})$ be a refinable function with a mask $m_0$. Define masks  $m_j: j=1,2,...,q$ so that  \\
1)$\hat\varphi (\chi\mathcal{A}^{-1})m_j(\chi)=
{\bf 1}_{E_j}(\chi)$, where  $E_j=G^\bot_{-s(j)}r_{-s(j)}^{\alpha_{-s(j)}}r_{-s(j)+1}^{\alpha_{-s(j)+1}}...r_{0}^{\alpha_{0}}...r_{M}^{\alpha_{M}}$ are disjoint cosets and  $E_j \mathcal{A}^t$ are disjoint also ,\\
2)there  are integers   $t(j)\ge 0$, such  that
$$
\bigsqcup_jE_j\mathcal{A}^{t(j)}=G_{M+1}^\bot\setminus G_M^\bot.
$$
 Then functions   $
\psi^{(1)},\psi^{(2)},...,\psi^{(q)}$ generate tight wavelet frame .
\end{theorem}

{\bf Proof.}

 We will calculate the sum

$$
\sum_{j=1}^q  \sum_{n=-\infty}^{+\infty}\sum_{h\in H_0}|(\psi^{(j)}_{n,h},f)|^2
$$
Using lemma  \ref{Lm3.1}
$$
\sum_{h\in H_0}|(\psi^{(j)}_{n,h},f)|^2=
\int_{E_j\mathcal{A}^{n}}|\hat{f}(\chi)|^2d\nu(\chi)
$$
we have. Therefore
$$
\sum_{n\in \mathbb Z}\sum_{h\in H_0}|(\psi^{(j)}_{n,h},f)|^2=
\sum_{n\in \mathbb Z}\int_{E_j\mathcal{A}^{n}}|\hat{f}(\chi)|^2d\nu(\chi)=
\sum_{n\in \mathbb Z}\int_{E_j\mathcal{A}^{n+t(j)}}|\hat{f}(\chi)|^2d\nu(\chi)
$$
Since sets  $E_j\mathcal{A}^{t(j)+n}$ - are disjoint ,
then

$$
\sum_{j=1}^q\sum_{n\in \mathbb Z}\sum_{h\in H_0}|(\psi^{(j)}_{n,h},f)|^2=
\sum_{n\in \mathbb Z}\sum_{j=1}^q\int_{E_j\mathcal{A}^{n+t(j)}}|\hat{f}(\chi)|^2d\nu(\chi)=
$$
$$
=\sum_{n\in \mathbb Z}\sum_{j=1}^q\int_{(E_j\mathcal{A}^{t(j)})\mathcal{A}^{n}}|\hat{f}(\chi)|^2d\nu(\chi)
=\sum_{n\in \mathbb Z}\int_{(\bigsqcup_{j}{E_j\mathcal{A}^{t(j)})\mathcal{A}^{n}}}|\hat{f}(\chi)|^2d\nu(\chi)
=$$
$$
=\sum_{n=-\infty}^{+\infty}\int_{(G_{M+1}^\bot\setminus G_M^\bot) \mathcal{A}^n }|\hat{f}(\chi)|^2d\nu(\chi)
=\int_X|\hat{f}(\chi)|^2d\nu(\chi).\square
$$
\subsection{Algorithm for constructing tight wavelet  frames. }
We will now indicate some class of  refinable step functions that generate tight wavelet  frames. Let $M=N \ge 1, p\ge 2$.\\
 {\bf First step.}
 In the $p$-adic tree $T$ (see Figure 2) we put
\begin{equation}\label{Eq3.2}
\lambda_{p^N}=\lambda_{p^N+1}=\dots=\lambda_{p^{N+1}-2} =0
\end{equation}
\begin{equation}\label{Eq3.3}
   \lambda_{p^{N+N+1}-1}=\lambda_{p^{N+N+1}-2}=\dots=\lambda_{p^{N+N+1}-p^N} =0.
   \end{equation}
   Denote the set of these numbers $\lambda_\nu$ as $\Lambda_0(T)$.
   We obtained  the initial tree  $T$ (see Figure 3).

   \unitlength=0.80mm
  \begin{picture}(240,80)
 \small

  \put(150,48){\line (-1,1){16}}
     \put(151,48){\line (0,1){16}}
     \put(152,48){\line (1,1){16}}
     \put(120,66){$\lambda_{p^{N+N+1}-p^N}=0$}
      \put(160,66){$\lambda_{p^{N+N+1}-1}=0$}

     \put(58,48){\line (-1,1){16}}
     \put(59,48){\line (0,1){16}}
     \put(60,48){\line (1,1){16}}
     \put(38,66){$\lambda_{p^{N+N}}$}
      \put(70,66){$\lambda_{p^{N+N}+p-1}$}
      \multiput(51,64)(3,0){6}{$\cdot$}

    \multiput(60,46)(3,0){31}{$\cdot$}
   \put(130,25){\line (-1,1){16}}
     \put(131,25){\line (0,1){16}}
     \put(132,25){\line (1,1){16}}
     \put(112,43){$0$}
     \put(130,43){$0$}
      \put(146,43){$\lambda_{p^{N+1}-1}$}
      \multiput(124,41)(3,0){6}{$\cdot$}
      \multiput(74,41)(3,0){4}{$\cdot$}

  \put(124,19){$\lambda_{p-1}$}
  \put(81,19){$\lambda_{1}$}
   \put(105,19){$\lambda_{\nu}$}
   \put(80,25){\line (-1,1){16}}
     \put(81,25){\line (0,1){16}}
     \put(82,25){\line (1,1){16}}
     \put(55,43){$\lambda_{p^{N}}=0$}
      \put(80,43){$0$}
       \put(96,43){$0$}
  \multiput(81,22)(3,0){17}{$\cdot$}

 \put(100,-2){$\lambda_{0}=1$}
  \put(100,2){\line (-1,1){16}}
  \put(106,2){\line (0,1){16}}
  \put(110,2){\line (1,1){16}}
  \put(-5,66){$\bf \mathfrak{G}_{N+1}^\bot \setminus \mathfrak{G}_{N}^\bot $:}
  \put(-5,43){$\bf \mathfrak{G}_{1}^\bot \setminus \mathfrak{G}_{0}^\bot $:}
  \put(-5,19){$\bf \mathfrak{G}_{-N+1}^\bot \setminus \mathfrak{G}_{-N}^\bot $:}
  \put(-5,-2){$\bf \mathfrak{G}_{-N}^\bot  $:}

    \end{picture}\\

\hskip4cm Figure 3. Initial  tree $T$.

\vskip0.4cm
   It is clear that $\sharp \Lambda_0(T)=p^{N+1}-1$.
 By theorem \ref{th2.01}: the  values $\lambda_\nu\in \Lambda_0(T)$ determine the mask $m_0$ of some  refinable  function $\varphi \in \mathfrak{D}_{\mathfrak G_{N}}(\mathfrak G_{-N})$.\\

   {\bf Second step.} We now define an elementary tree transformation as follows:\\
(i) Select a subtree

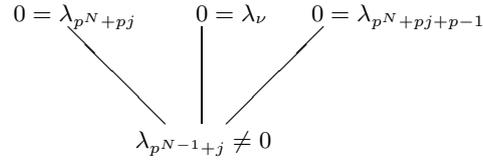
\unitlength=0.80mm
  \begin{picture}(240,30)
 \small
\put(95,-2){$\lambda_{p^{N-1}+j}\neq 0$}
  \put(100,2){\line (-1,1){16}}
  \put(106,2){\line (0,1){16}}
  \put(110,2){\line (1,1){16}}
  \put(124,19){$0=\lambda_{p^N+pj+p-1}$}
  \put(75,19){$0=\lambda_{p^N+pj}$}
   \put(105,19){$0=\lambda_{\nu}$}
   \end{picture}\\

and replace  it to a subtree

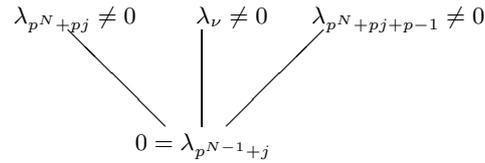
\unitlength=0.80mm
  \begin{picture}(240,30)
 \small
\put(95,-2){$0=\lambda_{p^{N-1}+j}$}
  \put(100,2){\line (-1,1){16}}
  \put(106,2){\line (0,1){16}}
  \put(110,2){\line (1,1){16}}
  \put(124,19){$\lambda_{p^N+pj+p-1}\neq 0$}
  \put(75,19){$\lambda_{p^N+pj}\neq 0$}
   \put(105,19){$\lambda_{\nu}\neq 0$}
   \end{picture}.\\

(ii) Select a subtree

\unitlength=0.80mm
  
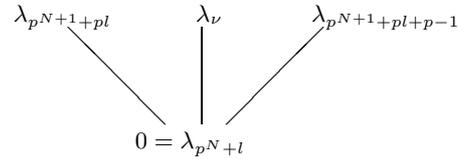
\begin{picture}(240,30)
 \small
\put(95,-2){$0=\lambda_{p^{N}+l}$}
  \put(100,2){\line (-1,1){16}}
  \put(106,2){\line (0,1){16}}
  \put(110,2){\line (1,1){16}}
  \put(124,19){$\lambda_{p^{N+1}+pl+p-1}$}
  \put(75,19){$\lambda_{p^{N+1}+pl}$}
   \put(105,19){$\lambda_{\nu}$}
   \end{picture}\\

and replace it to a subtree

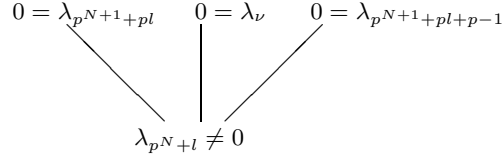
\unitlength=0.80mm
  \begin{picture}(240,30)
 \small
\put(95,-2){$\lambda_{p^{N}+l}\neq 0$}
  \put(100,2){\line (-1,1){16}}
  \put(106,2){\line (0,1){16}}
  \put(110,2){\line (1,1){16}}
  \put(124,19){$0=\lambda_{p^{N+1}+pl+p-1}$}
  \put(75,19){$0=\lambda_{p^{N+1}+pl}$}
   \put(105,19){$0=\lambda_{\nu}$}
   \end{picture}.\\

As a result we get   the new tree  $T$ for which $\sharp \Lambda_0(T)=p^{N+1}-1$.
Repeating this elementary transformation, we will obtain various trees for which $\sharp \Lambda_0(T)=p^{N+1}-1$.\\

{\bf Third  step.} Let  $T$ be a $p$-adic tree  received after the 2-nd step.  We transform the constructed tree as follows. For each node $\lambda_m\quad (m=p^n+\alpha_0p^0+\alpha_1p^1+...+\alpha_np^n, \alpha_n\neq p-1 )$, choose the path to the root
$$
\lambda_m\rightarrow \,\lambda_{m \,{\rm div}\, p} \rightarrow \lambda_{m \,{\rm div}\, p^2}\rightarrow\dots \rightarrow  \lambda_0=1.
$$
and replace $\lambda_m$ to the product
$$
\lambda_m\lambda_{m \,{\rm div}\, p}\lambda_{m \,{\rm div}\, p^2}\dots   \lambda_0 =\hat{\varphi}(\mathfrak{G}_{-N}^\bot r_{-N}^{\alpha_0}r_{-N+1}^{\alpha_1}...r_{-N+n}^{\alpha_n} )=:\hat{\varphi}_m,\quad \alpha_n\neq p-1.
$$

We have obtained a tree $T(\hat{\varphi})$ of values of the Fourier transform $\hat{\varphi}$.

\unitlength=0.80mm
  \begin{picture}(240,80)
 \small
  \put(150,48){\line (-1,1){16}}
     \put(151,48){\line (0,1){16}}
     \put(152,48){\line (1,1){16}}
     \put(125,66){$\hat{\varphi}_{p^{N+N+1}-p}$}
      \put(160,66){$\hat{\varphi}_{p^{N+N+1}-1}$}

     \put(58,48){\line (-1,1){16}}
     \put(59,48){\line (0,1){16}}
     \put(60,48){\line (1,1){16}}
     \put(38,66){$\hat{\varphi}_{p^{N+N}}$}
      \put(70,66){$\hat{\varphi}_{p^{N+N}+p-1}$}
      \multiput(51,64)(3,0){6}{$\cdot$}

    \multiput(60,46)(3,0){31}{$\cdot$}
   \put(130,25){\line (-1,1){16}}
     \put(131,25){\line (0,1){16}}
     \put(132,25){\line (1,1){16}}
     \put(110,43){$\hat{\varphi}_{p^{N}-p}$}
      \put(146,43){$\hat{\varphi}_{p^{N}-1}$}
      \multiput(124,41)(3,0){6}{$\cdot$}
      \multiput(74,41)(3,0){4}{$\cdot$}

  \put(124,19){$\hat{\varphi}_{p-1}$}
  \put(81,19){$\hat{\varphi}_{1}$}
   \put(105,19){$\hat{\varphi}_{\nu}$}
   \put(80,25){\line (-1,1){16}}
     \put(81,25){\line (0,1){16}}
     \put(82,25){\line (1,1){16}}
     \put(60,43){$\hat{\varphi}_{p^{N-1}}$}
      \put(86,43){$\hat{\varphi}_{p^{N-1}+p-1}$}

  \multiput(81,22)(3,0){17}{$\cdot$}

 \put(100,-2){$\hat{\varphi}_{0}=1$}
  \put(100,2){\line (-1,1){16}}
  \put(106,2){\line (0,1){16}}
  \put(110,2){\line (1,1){16}}
  \put(-5,66){$\bf \mathfrak{G}_{N+1}^\bot \setminus \mathfrak{G}_{N}^\bot $:}
  \put(-5,43){$\bf \mathfrak{G}_{0}^\bot \setminus \mathfrak{G}_{-1}^\bot $:}
  \put(-5,19){$\bf \mathfrak{G}_{-N+1}^\bot \setminus \mathfrak{G}_{-N}^\bot $:}
  \put(-5,-2){$\bf \mathfrak{G}_{-N}^\bot  $:}

    \end{picture}\\

\hskip4cm Figure 4.  The graph of the tree $T(\hat{\varphi})$.

\vskip0.4cm
{\bf Foutrh  step.}
 We transform the  tree $T(\hat{\varphi})$ as follows. In each subtree

 \unitlength=0.80mm
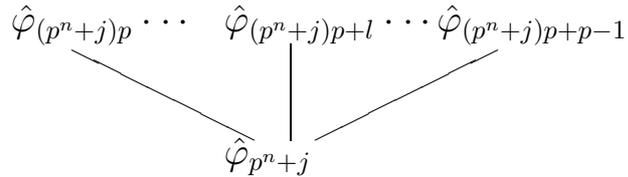
  \begin{picture}(240,30)
\put(95,-2){$\hat{\varphi}_{p^n+j}$}
  \put(100,2){\line (-2,1){30}}
  \put(106,2){\line (0,1){16}}
  \put(110,2){\line (2,1){30}}
  \put(130,20){$\hat{\varphi}_{(p^n+j)p+p-1}$}
  \put(60,20){$\hat{\varphi}_{(p^n+j)p } \cdot  \cdot \, \cdot$}
   \put(95,20){$\hat{\varphi}_{(p^n+j)p+l}\cdot  \cdot \, \cdot$}
   \end{picture}
   \vskip0.4cm
    \noindent  with $n=2N-1,2N-2,...,0$
 we replace $\hat{\varphi}_{(p^n+j)p+l},\quad (l=0,...,p-1)$ to $\hat{\varphi}_{p^n+j}$. We replace also $\hat{\varphi}_1, \hat{\varphi}_2,...,\hat{\varphi}_{p-1}$ to $\hat{\varphi}_0=1$.

  Thus, we have obtained a tree  of values of the Fourier transform $\hat{\varphi}(\chi\mathcal{A}^{-1})$.  Now we can specify the algorithm for constructing masks.
  There is a smallest level in the tree $T(\hat{\varphi}(\chi\mathcal{A}^{-1}))$ that contains at least one zero. Let it be the $n+1$-th level.
  Then the tree $T(\hat{\varphi}(\chi\mathcal{A}^{-1}))$ does not contain zeros in the $n$-th level.
  We know that the function $\hat{\varphi}(\chi\mathcal{A}^{-1})$ is constant on cosets with respect to the subgroup $\mathfrak{G}_{-N+1}^\bot$.\\
   If $n>1$, we represent  the set
   $\mathfrak{G}_{-N+n}^\bot \setminus \mathfrak{G}_{-N+n-1}^\bot $
   as the union of disjoint cosets
   $$
     \mathfrak{G}_{-N+m}^\bot r_{-N+m}^{\alpha_{-N+m}}r_{-N+m+1}^{\alpha_{-N+m+1}}...r_{-N+s}^{\alpha_{-N+s}}
     \mathcal{A}^{\gamma_{m,s}},\, s\ge m\ge 0.
      $$
      Then
      $$
      \bigsqcup\limits_{m,s} \mathfrak{G}_{-N+m}^\bot r_{-N+m}^{\alpha_{-N+m}}r_{-N+m+1}^{\alpha_{-N+m+1}}...r_{-N+s}^{\alpha_{-N+s}}
     \mathcal{A}^{\gamma_{m,s}}\mathcal{A}^{2N-n+1}=
     \mathfrak{G}_{N+1}^\bot \setminus \mathfrak{G}_{N}^\bot.
      $$
      Now we define functions $m_j(\chi)\hat{\varphi}(\chi\mathcal{A}^{-1})$ as characteristic functions of cosets
      \begin{equation}\label{Eq3.4}
      \mathfrak{G}_{-N+m}^\bot r_{-N+m}^{\alpha_{-N+m}}r_{-N+m+1}^{\alpha_{-N+m+1}}...r_{-N+s}^{\alpha_{-N+s}}.
           \end{equation}
             This means that functions
      $\hat{\psi}^{(j)}(\chi)$ are the characteristic functions of the cosets (\ref{Eq3.4} ) and functions $\psi^{(j)}  $ generate tight wavelet frame. \\
      If $n=1$ and $\lambda_1=\lambda_2=...=\lambda_{p-1}=0$ then
 $$
      m_0(\chi)={\bf 1}_{\mathfrak{G}_{-N}^\bot}(\chi),\,
      \hat{\varphi}(\chi)={\bf 1}_{\mathfrak{G}_{-N}^\bot}(\chi),\,
       \varphi(x)=p^{-N}{\bf 1}_{\mathfrak{G}_{-N}}(x)
 $$
 and function $\phi(x)=p^N\varphi(\mathcal{A}^N x)$ generates an orthogonal MRA.\\
 If $n=1$ and not all numbers  $\lambda_1, \lambda_2,...,\lambda_{p-1}$ are equal to zero then we define sets
 $$
 J_1=\{j=\overline{1,p-1}:\lambda_j\neq 0\}, \quad J_0=\{j=\overline{1,p-1}:\lambda_j = 0.\}
 $$
 Let $J_1=J_{1,-N}\bigsqcup J_{1,-N+1}$ be a disjoint partition of the set $J_1$.
       For   $j\in J_{1,-N+1}$ we define masks $m_j$ and wavelets $\psi^{(j)}$ by the equalities
       $$m_j(\chi)\hat{\varphi}(\chi\mathcal{A}^{-1})={\bf 1}_{\mathfrak{G}_{-N+1}^\bot r_{-N+1}^j}(\chi) =\hat{\psi}^{(j)}(\chi).
       $$
       For $j\in J_{1,-N}$ we define masks $m_{j,k}$ and wavelets $\psi^{(j,k)},\ k=\overline{0,p-1}$ by the equalities
       $$m_{j,k}(\chi)\hat{\varphi}(\chi\mathcal{A}^{-1})={\bf 1}_{\mathfrak{G}_{-N}^\bot r_{-N}^k r_{-N+1}^j}(\chi) =\hat{\psi}^{(j,k)}(\chi).$$
       For   $j\in J_0$ we define masks $m_j$ and wavelets $\psi^{(j)}$ by the equalities
       $$m_j(\chi)\hat{\varphi}(\chi\mathcal{A}^{-1})={\bf 1}_{\mathfrak{G}_{-N}^\bot r_{-N}^j}(\chi)=\hat{\psi}^{(j)}(\chi).
       $$
Since

\begin{multline}
\mathfrak{G}_{-N+1}^\bot\setminus\mathfrak{G}_{-N}^\bot=\\\biggl(\bigsqcup_{j\in J_{1,-N}}\mathfrak{G}_{-N}^\bot r_{-N}^k r_{-N+1}^j\biggr)\bigsqcup
\biggl(\bigsqcup_{j\in J_{1,-N+1}}\mathfrak{G}_{-N+1}^\bot r_{-N+1}^j\biggr)\bigsqcup
\biggl(\bigsqcup_{j\in J_0}\bigl(\mathfrak{G}_{-N}^\bot r_{-N}^j\bigr)\mathcal{A}\biggr)
\end{multline}
wavelets $\psi^{(j)}$ and $\psi^{(j,k)}$ generate a tight wavelet frame.

{\bf Example.} Let $p=3, M=N=1$.
Let's build an initial tree  $T$ (see Figure 3) and perform an elementary transformation. We obtain  the new tree $T(m_0)$. Using this tree $T(m_0)$, we build a Fourier transform of refinable  function $\varphi$ (see Figures 5,6).

\noindent
\unitlength=0.6mm
 \begin{picture}(300,45)
 \put(0,20){\line(1,0){270}}
 \qbezier(0,20)(15,15)(30,20)
 \qbezier(0,20)(5,15)(10,20)
\multiput(0,18)(30,0){9}{\line(0,1){4}}
   \put(0,30){\line(1,0){10}}
   \put(20,30){\line(1,0){10}}
   \put(10,23){\line(1,0){10}}
     \put(5,33){\normalsize 1}
   \put(15,25){\normalsize 0}
      \put(21,33){\normalsize $\neq 0$}
   \put(45,26){\normalsize 0}
\put(90,23){\line(1,0){180}}
\put(30,23){\line(1,0){40}}
\put(70,30){\line(1,0){20}}
\put(210,26){\normalsize 0}
\put(78,32){\normalsize $\neq 0$}
 \qbezier(0,20)(45,10)(90,20)
 \put(42,9){\small ${\mathfrak{G}_{1}^\bot}$}
 \put(12,9){\small ${\mathfrak{G}_{0}^\bot}$}
 \put(0,9){\small ${\mathfrak{G}_{-1}^\bot}$}
 \qbezier(90,20)(135,10)(180,20)
 \put(125,9){\small ${\mathfrak{G}_{1}^\bot}r_1$}
 \qbezier(180,20)(225,10)(270,20)
 \put(205,9){\small ${\mathfrak{G}_{1}^\bot}r_1^2$}
 \qbezier(0,20)(135,-10)(270,20)
 \put(130,0){\small $\mathfrak{G}_2^\bot $}
 \put(250,8){$\normalsize \hat{\varphi}(\chi)$}
 \end{picture}

\hskip2cm Figure 5. Fourier transform {$\normalsize \hat{\varphi}(\chi)$}

\noindent
\unitlength=0.6mm
 \begin{picture}(300,45)
 \put(0,20){\line(1,0){270}}
 \qbezier(0,20)(15,15)(30,20)
 \qbezier(0,20)(5,15)(10,20)
\multiput(0,18)(30,0){9}{\line(0,1){4}}
   \put(0,30){\line(1,0){30}}
     \put(15,33){\normalsize 1}
   \put(45,26){\normalsize 0}
\put(90,23){\line(1,0){120}}
\put(210,30){\line(1,0){60}}
 \put(240,32){\normalsize $\neq 0$}
\put(30,23){\line(1,0){30}}
\put(60,30){\line(1,0){30}}
\put(160,26){\normalsize 0}
\put(78,32){\normalsize $\neq 0$}
 \qbezier(0,20)(45,10)(90,20)
 \put(42,9){\small ${\mathfrak{G}_{1}^\bot}$}
 \put(12,9){\small ${\mathfrak{G}_{0}^\bot}$}
 \put(0,9){\small ${\mathfrak{G}_{-1}^\bot}$}
 \qbezier(90,20)(135,10)(180,20)
 \put(125,9){\small ${\mathfrak{G}_{1}^\bot}r_1$}
 \qbezier(180,20)(225,10)(270,20)
 \put(205,9){\small ${\mathfrak{G}_{1}^\bot}r_1^2$}
 \qbezier(0,20)(135,-10)(270,20)
 \put(130,0){\small $\mathfrak{G}_2^\bot $}
 \put(250,8){$\normalsize \hat{\varphi}(\chi \mathcal{A}^{-1})$}
 \end{picture}

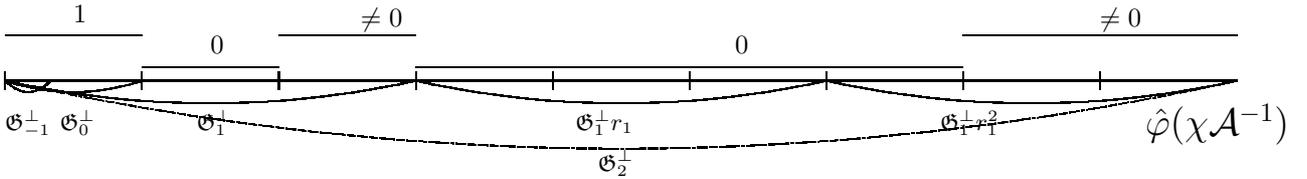
\hskip2cm Figure 6. Fourier transform {$\normalsize \hat{\varphi}(\chi\mathcal{A}^{-1})$}
\vskip0.5cm
Now we can write out three wavelet systems generating a tight wavelet frame:\\
1)
$$\hat{\psi}^{(1)}(\chi)=m_1(\chi)\hat{\varphi}(\chi\mathcal{A}^{-1})={\bf 1}_{\mathfrak{G}_{-1}^\bot r_{-1}^j}(\chi),
\hat{\psi}^{(2)}(\chi)={\bf 1}_{\mathfrak{G}_{0}^\bot r_{0}^2}(\chi),
$$
2)
$$\hat{\psi}^{(1)}(\chi)=
{\bf 1}_{\mathfrak{G}_{-1}^\bot r_{-1}}(\chi),\,
\hat{\psi}^{(2)}(\chi)={\bf 1}_{\mathfrak{G}_{1}^\bot r_{0}}(\chi),
$$
$$\hat{\psi}^{(3)}(\chi)=
{\bf 1}_{\mathfrak{G}_{-1}^\bot r_{-1}r_0}(\chi),\,
\hat{\psi}^{(4)}(\chi)={\bf 1}_{\mathfrak{G}_{1}^\bot r_{-1}^2r_{0}}(\chi),
$$
3)
$$\hat{\psi}^{(1)}(\chi)=
{\bf 1}_{\mathfrak{G}_{-1}^\bot r_{0}^2}(\chi),\,
\hat{\psi}^{(2)}(\chi)={\bf 1}_{\mathfrak{G}_{0}^\bot r_{0}r_{1}^2}(\chi),\,
$$
$$
\hat{\psi}^{(3)}(\chi)=
{\bf 1}_{\mathfrak{G}_{0}^\bot r_{1}^2r_0^2}(\chi),\,
\hat{\psi}^{(4)}(\chi)={\bf 1}_{\mathfrak{G}_{-1}^\bot r_{-1}}(\chi).
$$
\vskip0.3cm
{\it Acknowledgements}: This work was supported by the Russian Science Foundation
 No  22-21-00037,   https://rscf.ru/project/22-21-00037.


 \end{document}